\documentclass[%
twocolumn
]{mpi2015-cscpreprint}

\usepackage{booktabs}
\usepackage{orcidlink}
\usepackage{fontawesome}
\usepackage[american]{babel}
\usepackage{tikz}
\usepackage{pgfplots}
\pgfplotsset{%
  yticklabel style={anchor=east, text width=2em, align=right},
  compat=1.16,
  every axis/.style={%
    axis background/.style={fill=white},
    xmajorgrids,
    ymajorgrids,
    xlabel near ticks,
    ylabel near ticks,
    axis on top = false,
    scale only axis = false
  },
  every axis plot/.append style={%
    line width=\linew,
  },
  every axis plot post/.append style={%
    every mark/.append style={line width=.1pt}
  }
}
\usepackage{graphicx}

\usepackage{amssymb}
\usepackage{amsthm}
\usepackage{amsmath}

\usepackage{caption}
\usepackage{subcaption}

\usepackage{lineno}
\usepackage{placeins}

\usepackage[normalem]{ulem}

\renewcommand{\vec}{\mathrm }
\usepackage[textsize=scriptsize]{todonotes}
\usepackage{cleveref}
\usepackage[software, hardware]{mymacros}
\usepackage{siunitx}

\graphicspath{{./ext_figures/}}

\renewcommand{\hat}[1]{{\widehat{#1}}}
\renewcommand{\tilde}[1]{{\widetilde{#1}}}

\renewcommand{\hx}{\widehat{\mathrm{x}} }

\newcommand{\linew}{1.5pt}

\begin{document}

\title{Application of operator inference to reduced-order modeling of constrained mechanical systems}

\author[$\ast$, $\dagger$]%
{Peter Benner \orcidlink{0000-0003-3362-4103}}
\author[$\dagger$,$\ast$, \faEnvelopeO]%
       {Yevgeniya Filanova~\orcidlink{0000-0002-8599-3747}}
\author[$\ast$]%
       {\authorcr Igor Pontes Duff~\orcidlink{0000-0001-6433-6142}}     
\author[$\ast$]%
{Jens Saak~\orcidlink{0000-0001-5567-9637}}

\affil[$\ast$]{%
         Max Planck Institute for Dynamics of Complex Technical Systems,
         Magdeburg, Germany}
     
\affil[$\dagger$]{Otto von Guericke University Magdeburg, Magdeburg, Germany}
\affil[\faEnvelopeO]{Corresponding author,\email{filanova@mpi-magdeburg.mpg.de}}
\abstract{Constrained mechanical systems occur in many applications, such as modeling of robots and other multibody systems. In this case, the motion is governed by a system of differential-algebraic equations (DAE), often with large and sparse system matrices. The problem dimension strongly influences the effectiveness of simulations for system analysis, optimization, and control, given limited computational resources. Therefore, we aim to obtain a simplified surrogate model with a few degrees of freedom that is able to accurately represent the motion and other important properties of the original high-dimensional DAE model.
Classical model reduction methods intrusively exploit the system matrices to construct the projection of the high-fidelity model onto a low-dimensional subspace. In practice, the dynamical equations are frequently an inaccessible part of proprietary software products.
In this work, we show an application of the non-intrusive operator inference (OpInf) method to DAE systems of index 2 and 3. Considering the fact that for proper DAEs there exists an ODE realization on the so-called hidden manifold, the OpInf optimization problem directly provides the underlying ODE representation of the given DAE system in the reduced subspace. A significant advantage is that only the DAE solution snapshots in a compressed form are required for identification of the reduced system matrices. Stability and interpretability of the reduced-order model is guaranteed by enforcing the symmetric positive definite structure of the system operators using semidefinite programming.
The numerical results demonstrate the implementation of the proposed methodology for different examples of constrained mechanical systems, tested for various loading conditions.}

\novelty{\begin{itemize}
\item Data-driven identification of the reduced second-order ODE representations of proper DAE systems, applied for mechanical models with position and velocity constraints;
\item Elimination of the algebraic conditions using displacement and velocity snapshot matrices, depending on the constraint type;
\item Selection of the number of snapshots based on the impulse response of the system.
	\end{itemize}
 }
\maketitle

\section{Introduction}
Transient analysis of vibrating structures is usually computationally
challenging due to the large number of elements in discretized
models. Therefore, model-order reduction (MOR) methods are often employed to replace the full-dimensional model with a smaller surrogate system~\cite{morAnt05, morBenF19}. Constraints on the movement of a system's components necessitate even
more sophisticated mathematical models, which involve adding algebraic equations for these
constraints and, consequently, turning the models into differential-algebraic equations (DAEs).

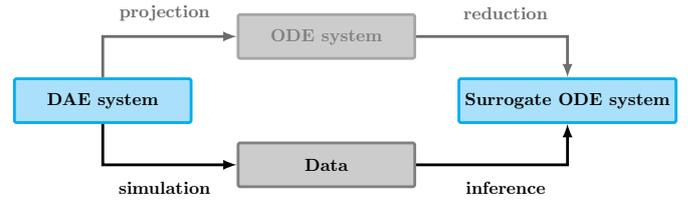
\begin{figure}[t]
  \centering
  \resizebox{0.49\textwidth}{!}{%
    \begin{tikzpicture}[ultra thick]
      \tikzstyle{every node}=[%
        font=\normalsize,
        minimum width = 100pt,
        minimum height = 25pt,
        rounded corners=.005\linewidth,
        font=\bfseries
      ]
      %
      %
      \node[fill=cyan!30!white, draw = cyan, minimum height = 25pt] (dae) {DAE system};
      %
      %
      \node [fill=gray!40!white, draw = gray!70!white, text = gray,%
             above right=10pt and 25pt of dae] (ode) {ODE system};
      %
      %
      \node [fill=gray!40!white, draw = gray, below right=10pt and 25pt of dae]
            (data) {Data};
      %
      %
      \node [fill=cyan!30!white, draw = cyan, right = 150pt of dae]
            (redode) {Surrogate ODE system};
      %
      %
      \draw[-latex, gray!85!black] (dae) |- node[align=center,above, xshift = 35pt] {projection} (ode);
      \draw[-latex] (dae) |- node[align=center,below, xshift = 35pt] {simulation} (data);
      \draw[-latex] (data) -| node[align=center,below, xshift = -35pt] {inference} (redode);
      \draw[-latex, gray!85!black] (ode) -| node[align=center,above, xshift = -35pt] {reduction} (redode);
    \end{tikzpicture}
  }%
  \caption{MOR strategies for proper DAE systems.}%
  \label{fig:scheme}
\end{figure}

\begin{figure*}[tp]
  \begin{subfigure}[b]{0.49\textwidth}
    \centering%
    \includegraphics[width=.9\linewidth]{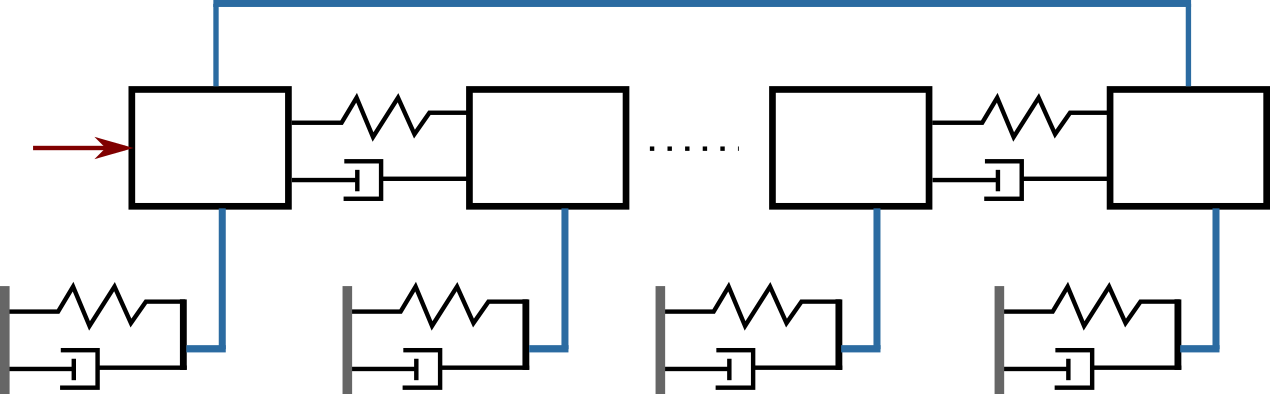}%
    \caption{Anchored mass-spring-damper system with connected masses~\cite{morMehS05}.}%
    \label{fig:examples_a}
  \end{subfigure}
  \hfill
  \begin{subfigure}[b]{0.49\textwidth}
    \centering%
    \includegraphics[width=.8\linewidth]{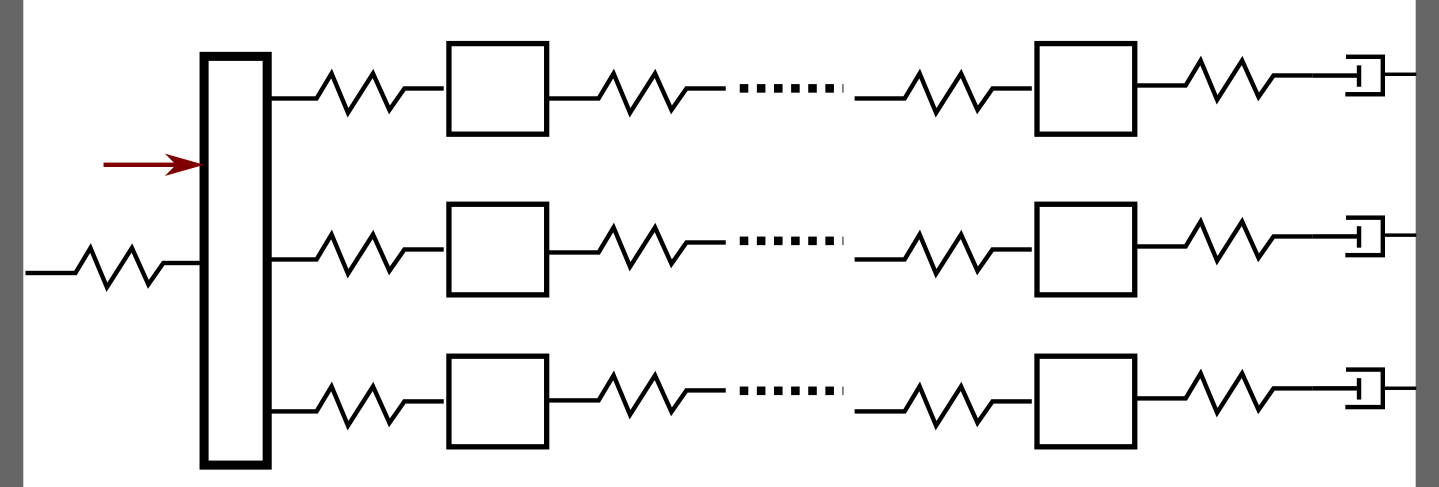}%
    \caption{Triple-chain mass-spring-damper system with velocity constraints~\cite{TruV09}.}%
    \label{fig:examples_b}
  \end{subfigure}
  \caption{Sketches of the two example systems.}%
  \label{fig:examples}
\end{figure*}

Many MOR methods are extended to DAE systems. Well-developed theories for
first-order descriptor systems exist, see~\cite{MehS06, morMehS05,
  gugercin2013model}. Constrained mechanical systems\footnote{cf. illustration in \Cref{fig:examples}} are typically modeled by second-order DAEs, like
\begin{align} \label{eq:mDAE}
  M \ddot{\vec{x}}(t) + D \dot{\vec{x}}(t) + K \vec{x}(t) + G_p^{\tran} \lambda_p(t) &+ G_v^{\tran} \lambda_v(t) = B \vec{u}(t), \nonumber \\
  G_p \vec{x}(t) &+ G_v \dot{\vec{x}} (t) = 0, \nonumber \\
  \vec{y}(t) = C_p \vec{x}(t) &+ C_v \dot{\vec{x}}(t),
\end{align}
which require structure-specific model order reduction approaches studied, e.g., in~\cite{morUdd15,morSaaV18,morAhmB14}. Often, the system~\eqref{eq:mDAE} can be represented by an equivalent ordinary differential equation system (ODEs) on a constraint manifold:
\begin{align}\label{eq:redODE}
  \hM \ddot{\hx} (t) + \hD \dot{\hx} (t) + \hK \hx (t) &= \hB \vec{u} (t), \\
  \vec{y} (t) &= \hC_v \dot{\hx} (t) + \hC_p \hx (t). \nonumber
\end{align}
In theory, the governing DAE system's matrices can be used to obtain the underlying ODE representation via projection, which can be further reduced using standard MOR techniques~\cite{morSaaV18}, as summarized in \Cref{fig:scheme}. However, even in the model-aware setting, the procedure requires a complicated analysis and adjustment of the existing solution algorithms to avoid the explicit calculation of projectors. Another drawback is the necessity
for calculations involving the system matrices, which can pose a significant challenge
when using commercial programs. Accessing the dynamical equations in their
(semi-)explicit form can be cumbersome or expensive also for open-source tools when
their discretization frameworks work matrix-free. Therefore, we focus on
alternative approaches, i.e., on identification of the reduced system from simulation data that is readily available in all systems. As shown in \Cref{fig:scheme}, we aim to obtain the reduced ODE system directly using the DAE solution data.

To this end, in this paper, we propose another application of the operator inference method~\cite{morPehW16} for a broad class of constrained mechanical systems that can be considered as proper index-2 and index-3 DAE systems. Due to the simplicity and flexibility of the operator inference method, its rapid development for various system structures also includes applications to DAEs. In~\cite{morKho22}, the authors describe an operator inference approach for DAE systems derived from lifting transformations of nonlinear models. The proposed approach enables the separate identification of the ODE part and the determination of algebraic operators. On the other hand, the work~\cite{morBenGHetal22} shows the identification of a reduced ODE system for semi-discretized Navier-Stokes equations. The algebraic constraints are eliminated through projection of the dataset onto a subspace spanned by the leading proper orthogonal decomposition (POD) modes.

We generalize the results presented in~\cite{morBenGHetal22} for a large group of mechanical systems with position or velocity constraints. The proposed approach allows learning the second-order ODE structure~\eqref{eq:redODE} directly from the time-domain data provided by a suitable DAE integrator.

The paper is organized as follows: \Cref{sec:mechsys} discusses the details of high-dimensional DAE models for constrained mechanical systems and the respective model-order reduction strategies; \Cref{sec:opinf} describes the main operator inference reduction procedure; The respective numerical results are presented in \Cref{sec:numres}.

\section{Model-order reduction of proper constrained mechanical systems}%
\label{sec:mechsys}
Multibody systems with constraints are often modeled by a second-order DAEs in the Lagrange multiplier form~\cite{eich1998numerical,Hau96,shabana2020dynamics,simeon2001numerical,rentrop99}, as shown in~\eqref{eq:mDAE}, where $M, D, K \in \mathbb{R}^{n \times n}$ are the symmetric positive definite mass, damping, and stiffness matrix, respectively. The state vector $\vec{x} \in \mathbb{R}^{n}$ represents displacements at each node of a space-discretized structure. Motion constraints are incorporated into the system via Lagrange multipliers represented by the vectors $\lambda_p, \; \lambda_v \in \mathbb{R}^{n_c}$, and full-rank Jacobian matrices $G_p$, $G_v$. External loading is applied through an input signal $\vec{u}(t) \in \mathbb{R}^{m}$ and a control matrix $B$. The system response is measured by the output vector $\vec{y}(t) \in \mathbb{R}^{p}$, which is related to the state and the first state derivative via the output matrices $C_p$ and $C_v$, respectively. In practice, it is common to have the constraints either only for the positions or only for the velocities. Hence, we will consider systems with only one of the matrices, $G_p$ or $G_v$.

The existence of the equivalent ODE representation for proper DAE systems, i.e., systems with a finite-limit transfer function, is proven by well-developed theories for first-order DAEs~\cite{morMehS05, kunkel2006differential}. Transforming the system~\eqref{eq:mDAE} into its first-order companion form doubles the system dimension and causes a loss of structure. However, simple observations allow us to verify that system~\eqref{eq:mDAE} possesses an underlying second-order ODE representation.

First, assume that the DAE system only has position constraints:
\begin{equation}
  G_p \vec{x}(t)= 0.
  \label{eq:pos_constr}
\end{equation}
Due to the full-rank condition on the constraint matrix $G_p$, restricting the system~\eqref{eq:mDAE} onto the kernel $\kernel{G_p}$, spanned by $N \in \mathbb{R}^{n \times n_{c}}$, leads to the system~\eqref{eq:proj_ker} with automatically satisfied algebraic restrictions~\eqref{eq:pos_constr}. With $\vec{x} (t) = N \vec{z} (t)$, we have
\begin{align}
  N^{\tran} M N \ddot{\vec{z}} (t) + N^{\tran} D N \dot{\vec{z}} (t) + N^{\tran} K N \vec{z} (t) = N^{\tran} B \vec{u} &(t), \nonumber \\
  \vec{y} (t) = C_p N \vec{z} (t) + C_v \dot{\vec{z}} &(t).
                \label{eq:proj_ker}
\end{align}
Analogously, in case of pure velocity constraints
\begin{equation}
  G_v \dot{\vec{x}}(t)= 0,
  \label{eq:vel_constr}
\end{equation}
the system~\eqref{eq:mDAE} can be projected onto $\kernel{G_v}$. A more detailed
description on obtaining the underlying ODE representations for second-order DAEs is given in~\cite{morSaaV18, morUdd15}. The explicit construction of the
projectors can be avoided by adjusting the standard reduction methods, as also
shown in~\cite{morSaaV18}, based on the derivations for first-order systems in~\cite{morHeiSS08, gugercin2013model}.

Another way to directly obtain the reduced ODE representation of~\eqref{eq:mDAE} is to apply the proper orthogonal decomposition (POD) method~\cite{morKunV08}. In this case, dominant characteristic structures are extracted from the solution trajectories of~\eqref{eq:mDAE}, collected at pre-defined time points $t_1, \dots , t_k$ in a snapshot matrix:
\begin{equation}
  X = \begin{bmatrix}
    | &  \cdots & | \\
    \vec{x}(t_1) & \cdots & \vec{x}(t_k) \\
    | & \cdots & |
  \end{bmatrix}.
  \label{mat:state_snap}
\end{equation}
We can find the low-dimensional subspace basis by performing a singular value decomposition (SVD) of the displacement snapshot matrix
\begin{equation}
X = V \Sigma W^{\tran}
\label{eq:svd_pos}
\end{equation}
and truncating the last $n - r$ left singular vectors in $V$. The resulting basis $V_r \in \mathbb{R}^{n \times r}$ spans a subspace of $\range{X}$. As discussed in~\cite{morBenGHetal22}, for the case of position constraints, i.e., $G_p \neq 0$ and $G_v = 0$,  it holds
\begin{equation}
G_p X = 0 \implies G_p V_r = 0.
\end{equation}
Analogously, for $G_p = 0$ and $G_v \neq 0$, we can define the reduced subspace basis $V_r$ by performing an SVD of the velocity snapshot matrix
\begin{equation}
\dot{X} = V_v \Sigma_v W_v^{\tran},
\label{eq:svd_vel}
\end{equation}
and truncating the left singular vectors $V_v$. The reduced basis $V_r$ also satisfies the algebraic conditions:
\begin{equation}
G_v \dot{X} = 0 \implies G_v V_r = 0.
\end{equation}

 Finally, projecting the system~\eqref{eq:mDAE} onto the respective POD subspace yields the reduced second-order ODE system:
\begin{align}\label{eq:redODE_POD}
\hM \ddot{\hx} (t) + \hD \dot{\hx} (t) + \hK \hx (t) &= \hB \vec{u} (t) \\
\vec{y} (t) &= \hC_v \dot{\hx} (t) + \hC_p \hx (t), \nonumber
\end{align} where the reduced state $\hx (t) \in \mathbb{R}^r, \; r \ll n$. The respective system matrices are
\begin{align}
	\hM = V^{\tran} M V,& \qquad \hD = V^{\tran} D V, \qquad \hK \, = V^{\tran} K V,   \nonumber \\
	 \hB = V^{\tran} B, \phantom{V} & \qquad \hC_v = C_v V, \phantom{V}  \qquad \hC_p = C_p V .
	\label{mat:pod}
\end{align}
The above described methods require access to the system matrices, which can be a key issue when working with proprietary software. Therefore, in this work, we propose to learn the reduced ODE representation~\eqref{eq:redODE_POD} from simulation data of the full-order DAE system~\eqref{eq:mDAE}, using the operator inference method~\cite{morPehW16}.

\section{Identification of the underlying reduced ODE system with the operator inference method}%
\label{sec:opinf}
As shown in \Cref{fig:scheme}, the reduced ODE system that approximates the dynamics of the given high-dimensional DAE system can be obtained by directly applying data-driven techniques. In the following, we use the second-order operator inference method described in~\cite{morFilPGetal23} extended for a control input-output system.

First, we assume to have access to the solution data, provided by an integrator that operates directly on the given DAE, e.g.,~\cite{Bot08, Hau96}, so that the computed numerical solution exactly satisfies the constraints. The time-domain data corresponding to the specific input signal $\vec{u} (t)$ is collected in the following snapshot matrices at pre-defined time points $t_1,\dots, t_k$:
\begin{align}
  \label{mat:snap}
  U &=
  \begin{bmatrix}
    | &  \cdots & | \\
    \vec{u}(t_1) & \cdots & \vec{u}(t_k) \\
    | & \cdots & |
  \end{bmatrix}, \quad
  Y =
  \begin{bmatrix}
    | &  \cdots & | \\
    \vec{y}(t_1) & \cdots & \vec{y}(t_k) \\
    | & \cdots & |
  \end{bmatrix}, \nonumber \\
  X &= \begin{bmatrix}
    | &  \cdots & | \\
    \vec{x}(t_1) & \cdots & \vec{x}(t_k) \\
    | & \cdots & |
  \end{bmatrix}.
\end{align}
Analogously, we assemble the respective derivatives of the state vector $\dot{X}, \ddot{X}$ that are available, when using~\cite{Bot08}.

Now, the reduced subspace can be identified from the collected data, as described in \Cref{sec:mechsys}. Note that when the integrator guarantees the fulfillment of the algebraic conditions, the POD modes $V_r$ also belong to a constrained manifold, defined by the algebraic equations. The compressed data set is calculated by restricting the snapshot matrices to the low-dimensional POD subspace:
\begin{equation} \label{eq:compressed_data}
X_r = V_r^{\tran} X, \quad \dot{X}_r = V_r^{\tran} \dot{X}, \quad \ddot{X}_r = V_r^{\tran} \ddot{X}.
\end{equation}

The final step is formulating the optimization problem and fitting the reduced operators to the compressed data according to~\eqref{eq:redODE_POD}. However, a good fitting of the operators, alone, does not guarantee a good dynamical system approximation. Therefore, it is necessary to enforce the properties inherent in the original system matrices, namely the symmetric positive (semi-)definite (SPD) structure:
\begin{equation}\label{eq:SPDconstr}
M \succ 0, \quad D \succeq 0, \quad K \succ 0.
\end{equation}
Structure preservation ensures stability, as well as interpretability of the results~\cite{morSalL06, Tho21}. Note, that only non-strict constraints can be enforced in the optimization problem. Therefore, the relations~\eqref{eq:SPDconstr} should be relaxed to symmetric positive semidefinite inequalities, as described in~\cite{morFilPGetal23}.

In this work, we treat the input and output matrices as unknown, unlike in force-informed operator inference, as described in~\cite{morFilPGetal23}. Formulating the optimization problem as a semidefinite program enables us to incorporate the SPD constraints~\eqref{eq:SPDconstr} and identify an optimal system representation that fits the given data.
Hence, the optimization process consists of two independent tasks. The state-related matrices are obtained by solving the following optimization problem:

\begin{align}\label{eq:sys_optim}
  \left[\hM ~~ \hD ~~ \hK ~~ \hB \right]  &= \nonumber \\ \underset{\substack{ M_r \succeq 0,  \; D_r \succeq 0, \\ K_r \succeq 0, \; B_r}}{\text{arg} \text{min}} &\left\| M_r \ddot{X}_r + D_r \dot{X}_r + K_r X_r - \tilde{B} U \right\|^2_\text{F}.
\end{align}
The output operators are identified analogously by solving a separate least-squares problem:
\begin{equation} \label{eq:y_optim}
\left[\hC_v \; \hC_p \right]  =  \underset{C_{v_r},\; C_{p_r}}{\text{arg} \text{min}} \left\| C_{v_r} \dot{X}_r + C_{p_r} X_r - Y \right\|^2_\text{F}.
\end{equation}
As a result, we learn the reduced ODE representation of the DAE system~\eqref{eq:mDAE} by solving the two minimization problems~\eqref{eq:sys_optim}
and~\eqref{eq:y_optim} using the reduced data~\eqref{eq:compressed_data}.

\section{Numerical results}\label{sec:numres}
\begin{figure*}[tbp]
  \begin{subfigure}[t]{.49\textwidth}
    \centering
    \includegraphics{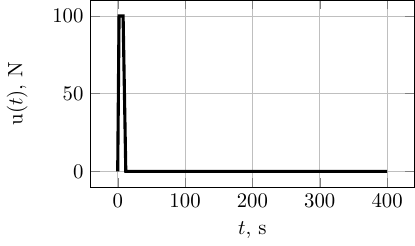}%
    \caption{Impulse-like input signal for snapshot generation.}%
    \label{fig:training_input}
  \end{subfigure}
  \begin{subfigure}[t]{.49\textwidth}
	\centering
	\includegraphics{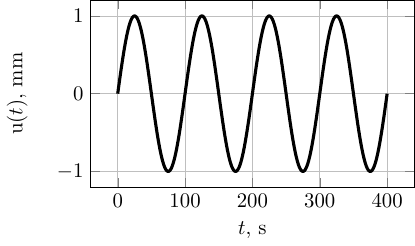}%
	\caption{Harmonic input signal for the ROM test.}%
	\label{fig:test_signal}
  \end{subfigure}
  \caption{Input signals for training and testing the ROMs.}%
  \label{fig:training_response}
\end{figure*}
\begin{figure*}[tp]
  \begin{subfigure}{0.49\textwidth}
  	\centering
    \includegraphics{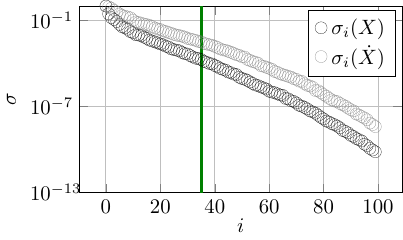}%
    \caption{Singular value decay.}%
    \label{fig:svd_ind2}
  \end{subfigure}
  \begin{subfigure}{0.49\textwidth}
  	\centering
    \includegraphics{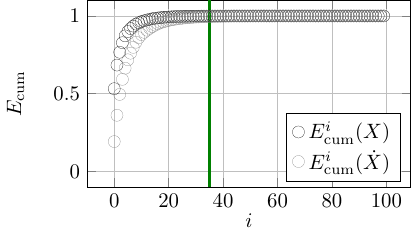}%
    \caption{Cumulative energy~\eqref{eq:cum_energy}.}%
    \label{fig:cum_ind2}
  \end{subfigure}
\centering
  \caption{Reduced order decision criteria for the triple-chain example with velocity constraints. The vertical line marks the selected reduced order \(r=35\).}%
  \label{fig:svd2}
\end{figure*}
\begin{figure*}[!ht]
	\begin{subfigure}{0.99\textwidth}
		\centering
		\includegraphics{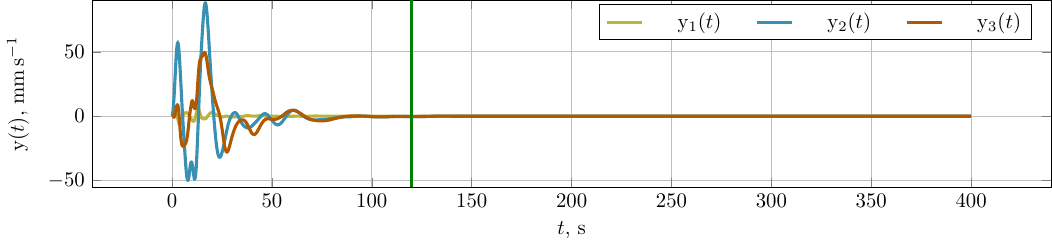}%
		\caption{Triple-chain example.}
	\end{subfigure}
	\vskip\baselineskip
	\begin{subfigure}{0.99\textwidth}
		\centering
		\includegraphics{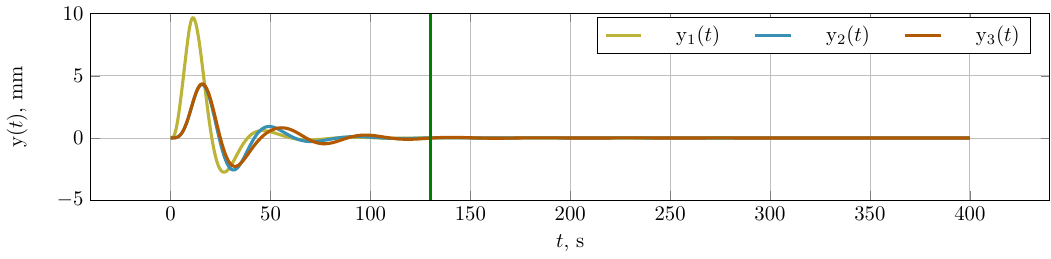}%
		\caption{Anchored-mass example.}
	\end{subfigure}
	\caption{System output response to the external load $\vec{u}(t)$, plotted in \Cref{fig:training_input}.}%
	\label{fig:training_response2}
\end{figure*}
\begin{figure*}[tp]
  \begin{subfigure}{0.49\textwidth}
  	\centering
    \includegraphics{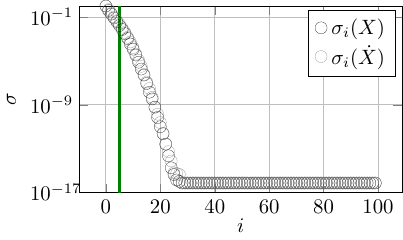}%
    \caption{Singular value decay.}%
    \label{fig:svd_ind3}
  \end{subfigure}
  \begin{subfigure}{0.49\textwidth}
  	\centering
    \includegraphics{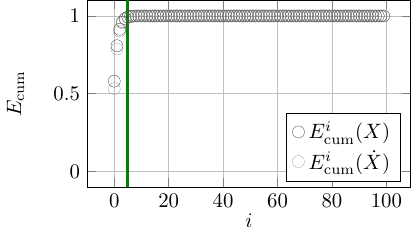}%
    \caption{Cumulative energy~\eqref{eq:cum_energy}.}%
    \label{fig:cum_ind3}
  \end{subfigure}
\centering
  \caption{Reduced order decision criteria for the anchored-mass example with position constraints.  The vertical line marks the selected reduced order \(r=5\).}%
  \label{fig:svd3}
\end{figure*}
\begin{figure*}[tp]
  \begin{subfigure}{1\textwidth}
  	\centering
    \includegraphics{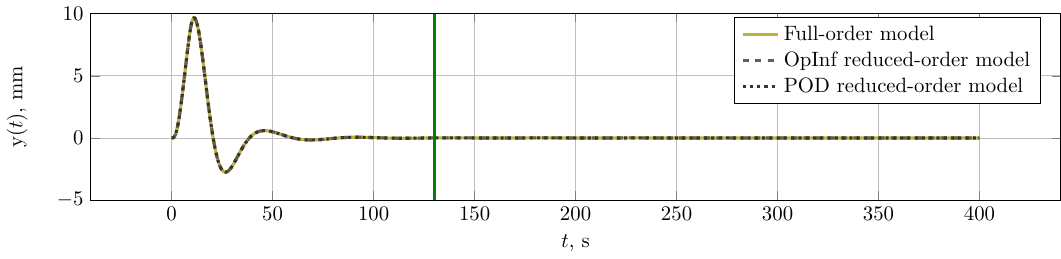}%
    \caption{First output.}
  \end{subfigure}
  \vskip\baselineskip%
  \begin{subfigure}{1\textwidth}
  	\centering
    \includegraphics{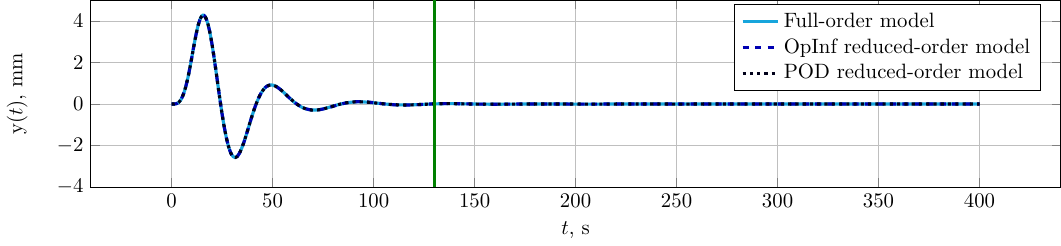}%
    \caption{Second output.}
  \end{subfigure}
  \vskip\baselineskip%
  \begin{subfigure}{1\textwidth}
  	\centering
    \includegraphics{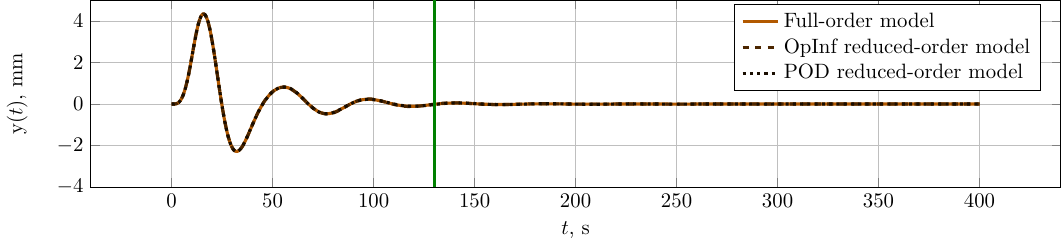}%
    \caption{Third output.}
  \end{subfigure}
  \vskip\baselineskip%
  \begin{subfigure}{\textwidth}
  	\centering
    \includegraphics{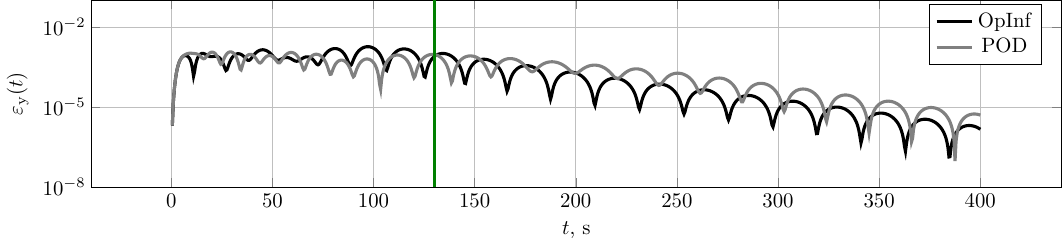}%
    \caption{Relative error of the ROM output trajectories.}%
    \label{fig:ind3_error_train}
  \end{subfigure}
  \caption{Validation of the ROM output trajectories for the anchored-mass example using the training input from \Cref{fig:training_input}. The vertical line marks the last snapshot, used for the ROM construction.}%
  \label{fig:ind3_train_signal_rom}
\end{figure*}
\begin{figure*}[tp]
  \begin{subfigure}{1\textwidth}
  	\centering
    \includegraphics{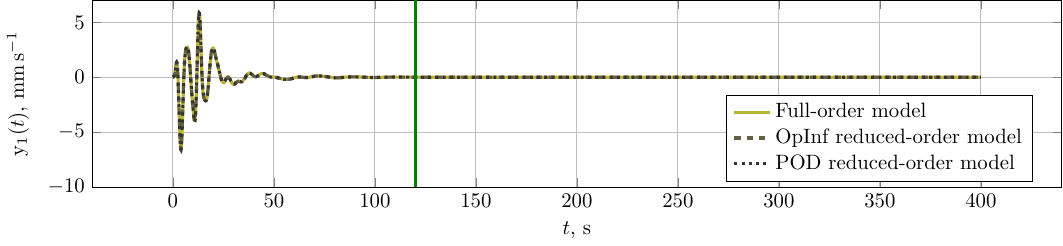}%
    \caption{First output.}
  \end{subfigure}
  \vskip\baselineskip%
  \begin{subfigure}{1\textwidth}
  	\centering
    \includegraphics{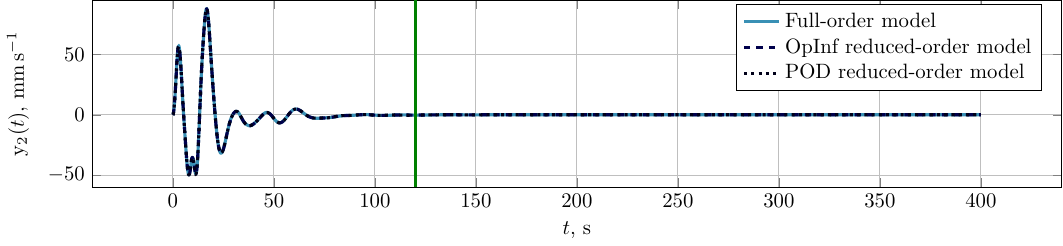}%
    \caption{Second output.}
  \end{subfigure}
  \vskip\baselineskip%
  \begin{subfigure}{1\textwidth}
  	\centering
    \includegraphics{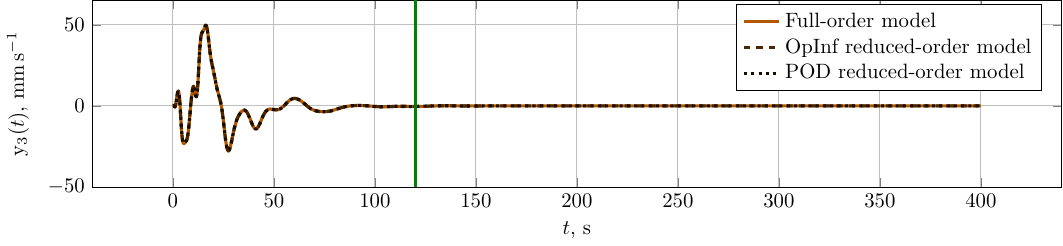}%
    \caption{Third output.}
  \end{subfigure}
  \vskip\baselineskip%
  \begin{subfigure}{\textwidth}
  	\centering
	\includegraphics{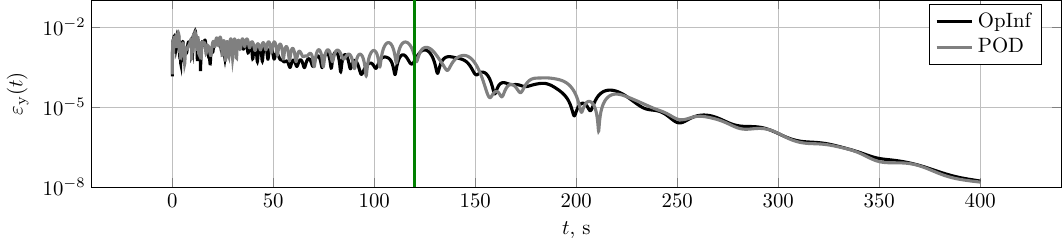}%
	\caption{Relative error of the ROM output trajectories.}%
	\label{fig:error_train}
  \end{subfigure}
  \caption{Validation of the ROM output trajectories for the triple-chain example using the training input from \Cref{fig:training_input}. The vertical line marks the last snapshot, used for the ROM construction.}%
  \label{fig:train_signal_rom}
\end{figure*}
\begin{figure*}[tp]
  \begin{subfigure}{1\textwidth}
  	\centering
    \includegraphics{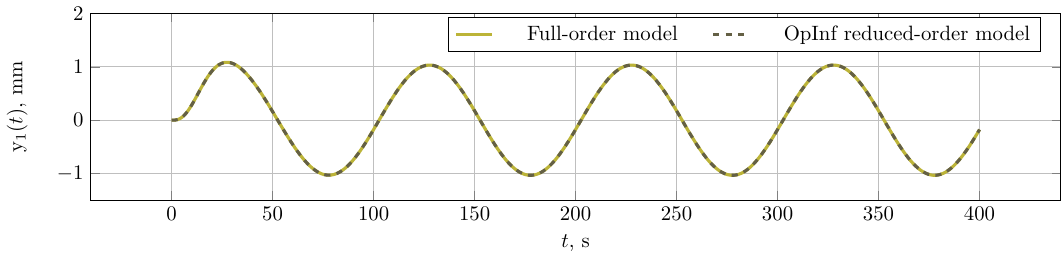}%
    \caption{First output.}
  \end{subfigure}
  \vskip\baselineskip%
  \begin{subfigure}{1\textwidth}
  	\centering
    \includegraphics{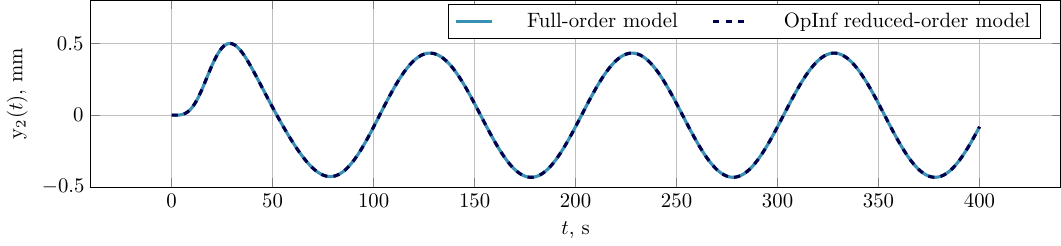}%
    \caption{Second output.}
  \end{subfigure}
  \vskip\baselineskip%
  \begin{subfigure}{1\textwidth}
  	\centering
    \includegraphics{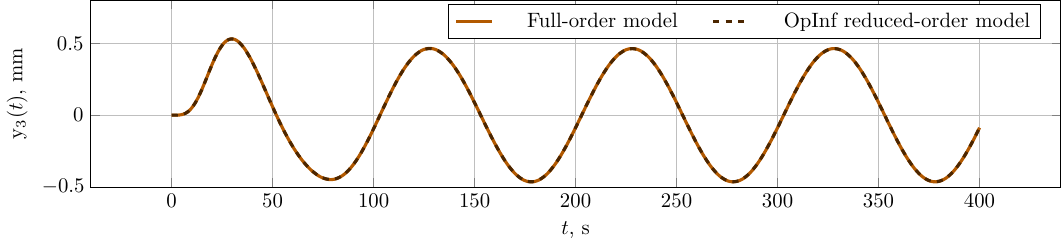}%
    \caption{Third output.}
  \end{subfigure}
  \vskip\baselineskip%
  \begin{subfigure}{\textwidth}
  	\centering
    \includegraphics{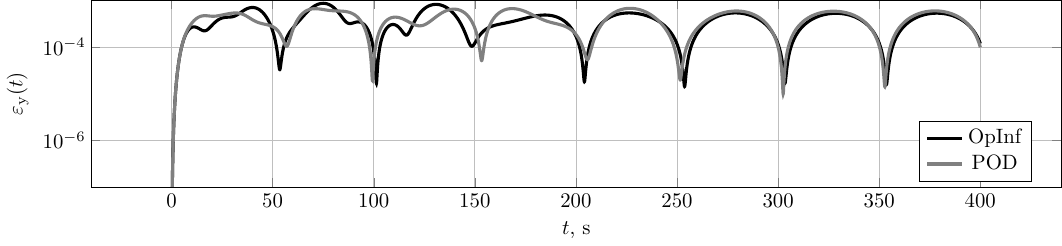}%
    \caption{Relative error for the test comparison of the ROM output trajectories.}%
    \label{fig:ind3_error_test}
  \end{subfigure}
  \caption{Comparison of the ROM output trajectories for the anchored-mass example using the harmonic test input from \Cref{fig:test_signal}.}%
  \label{fig:ind3_test_signal_rom}
\end{figure*}
\begin{figure*}[tp]
  \begin{subfigure}{1\textwidth}
  	\centering
    \includegraphics{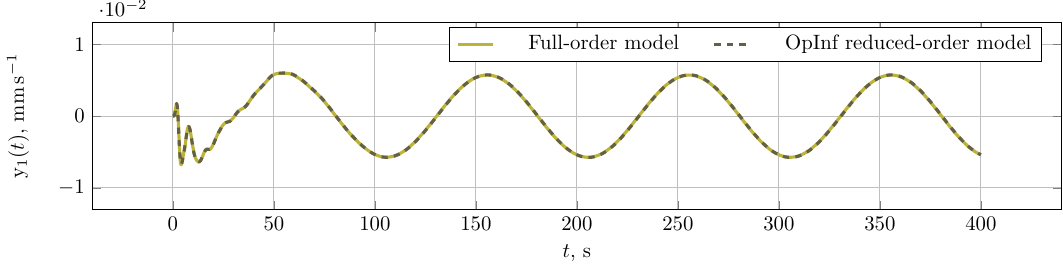}%
    \caption{First output.}
  \end{subfigure}
  \vskip\baselineskip%
  \begin{subfigure}{1\textwidth}
  	\centering
    \includegraphics{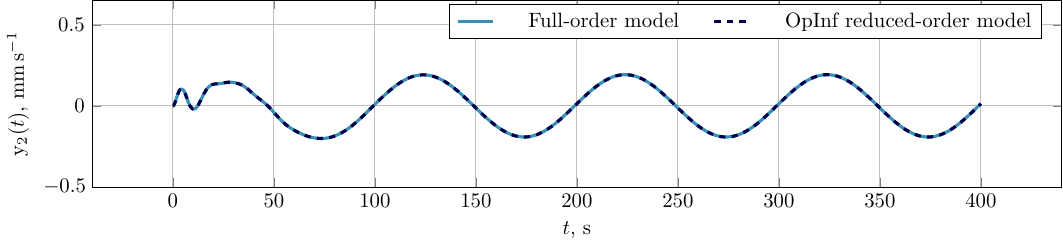}%
    \caption{Second output.}
  \end{subfigure}
  \vskip\baselineskip%
  \begin{subfigure}{1\textwidth}
  	\centering
    \includegraphics{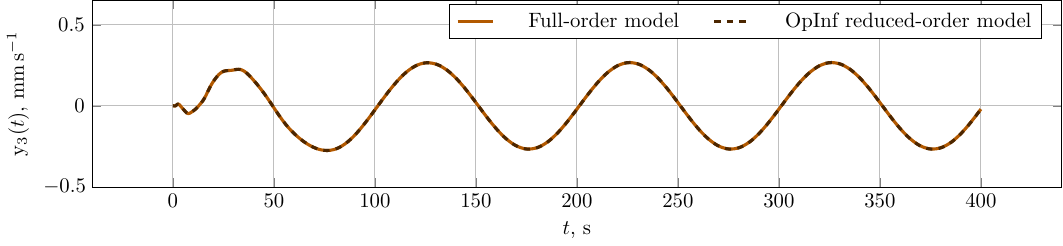}%
    \caption{Third output.}
  \end{subfigure}
  \vskip\baselineskip%
  \begin{subfigure}{\textwidth}
  	\centering
    \includegraphics{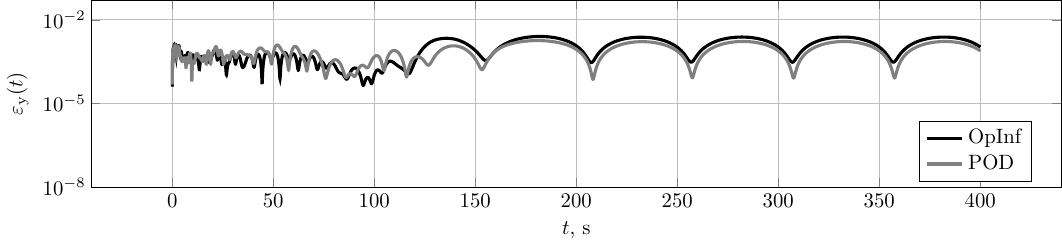}%
    \caption{Relative error for the test comparison of the ROM output trajectories.}%
    \label{fig:error_test}%
  \end{subfigure}
  \caption{Test comparison of the ROM output trajectories for the triple-chain example using the harmonic test input from \Cref{fig:test_signal}.}%
  \label{fig:test_signal_rom}
\end{figure*}
\subsection{Data availability and system requirements}
The simulation data, and code to produce the results presented in this paper are openly available at
 \url{https://doi.org/10.5281/zenodo.15724263}.
 The numerical experiments were run on a Dell computer with a 12th Gen \intel\coreifive-12600K (6 performance cores, 4 economy cores) and 32 GB RAM\@. The system was running Ubuntu version 20.04.4 (64Bit). The procedure described in \Cref{sec:opinf} is implemented in Matlab R2021a, using the YALMIP Toolbox R20210331~\cite{Lofberg2004}.
\subsection{Examples}
 In the following, we consider two mass-spring-damper (MSD) oscillators. The first model is shown in \Cref{fig:examples_a} and corresponds to an index-3 DAE system due to the imposed position constraints~\cite{morMehS05}. The masses in this example are anchored to the ground, the endpoints are additionally connected by a rigid bar. The second model in \Cref{fig:examples_b} consists of three rows of masses with velocity constraints (triple-chain example)~\cite{TruV09}.
The DAEs are integrated with a second-order Newmark integrator with embedded algebraic constraints~\cite{Bot08}. The reduced ODEs are integrated with the classical Newmark method~\cite{New59}.

According to the procedure in \Cref{sec:opinf}, we need to define
\begin{enumerate}
\item the simulation parameters for the snapshot set, such as external input function $\vec{u}(t)$, time step $dt$, end time $t_k$, etc.,
\item reduced order $r$,
\item analysis parameters, such as values for comparison, error measure $\varepsilon$, etc.
\end{enumerate}

The simulation parameters should be selected in a way to produce a data set that is informative but not redundant. We propose exciting the system with an impulse-like input signal and defining the step size and end time based on the system's output response. Thus, both oscillators are loaded with $\vec{u}(t)$, as illustrated in \Cref{fig:training_input}. Analyzing the corresponding output behavior allows us to select $dt$ and $t_k$ such that the trajectories are accurately captured and significant vibrations are faded out. In \Cref{fig:training_response2}, the corresponding output trajectories are shown. The vertical line represents the selected end time point $t_k$ for the snapshot generation.

Now, having the data set, we form the snapshots matrices, as shown in~\eqref{mat:snap}. The singular value decay of the snapshot matrix helps to understand the underlying dynamics and to choose the appropriate reduced order. For a more informed decision, we calculate the cumulative energy, as follows:
\begin{equation}
E^i_{\text{cum}} = \frac{\sum_{n=1}^{i}\sigma_n}{\sum_{n=1}^{N} \sigma_n}.
\label{eq:cum_energy}
\end{equation} 
We analyze the singular value decay of the state and derivative snapshot matrices according to \eqref{eq:svd_pos} and \eqref{eq:svd_vel}, which is particularly important for the MSD system with velocity constraints (see \Cref{fig:svd_ind2} and \Cref{fig:cum_ind2}). To construct the reduced basis $V_r$, we truncate the left singular vectors once $E_{\text{cum}}$ reaches its saturation point, which is shown in \Cref{fig:cum_ind2} and \Cref{fig:cum_ind3}. Hence, for the index-3 DAE system with position constraints of dimension $n = 600$, we select $r = 5$; for the index-2 DAE system with velocity constraints of dimension $n = 301$ we pick $r = 35$.

First, we compare the output trajectories of the reduced-order models (ROMs) and the full-order model (FOM) corresponding to the training input signal (see \Cref{fig:training_input}), as well as to the test input signal, a simple harmonic wave (see \Cref{fig:test_signal}).
 The comparison is performed by estimating the relative error in the output response for each time point $t_i$:
\begin{align}\label{eq:relerror}
\varepsilon_{\text{y}} (t_i) = \frac{ \| \vec{y} (t_i) - \widehat{\vec{y}}(t_i) \|_2^2 }{\underset{t \in [t_1, t_k]}{\max} \| \vec{y} (t) \|_2^2}.
\end{align} This type of error measurement enables adequate comparison of output trajectories around the zero value.

The results are presented according to the parameters and validation strategy described above. The comparison of the outputs according to the training input signal is shown in \Cref{fig:train_signal_rom} and \Cref{fig:ind3_train_signal_rom} for the MSD system with velocity and position constraints, respectively.
 We observe an accurate representation of the position output trajectories for the MSD system with position constraints based on the error plot in \Cref{fig:ind3_error_train}. The system dynamics can be captured by using only $r = 5$ POD modes, which clearly follows from the rapid singular value decay in \Cref{fig:svd_ind3}.

The MSD model with velocity constraints exhibits a slower singular value decay, as shown in \Cref{fig:svd_ind2}, which leads to the ROM of dimension $r = 35$. The more unknown variables are introduced for the optimization problem~\eqref{eq:sys_optim}, the more difficult it is to solve and get an accurate optimum. However, in \Cref{fig:error_train} the error value remains small and does not increase further after the vibrations have almost faded out. We can also see that the minor fluctuations are approximated with less accuracy than the significant vibrations at the start.

Finally, in \Cref{fig:ind3_test_signal_rom} and \Cref{fig:test_signal_rom}, we demonstrate the performance of the identified ROMs for the test input signal from \Cref{fig:test_signal}. The plots in \Cref{fig:ind3_error_test} and \Cref{fig:error_test} show that the constructed ROMs accurately approximate the dynamic behavior of the original DAE model when the external load differs from the training load used for the learned snapshot set.

\section{Conclusion}
In this work, we proposed an application of the operator inference methodology adapted for mechanical systems with position and velocity constraints. The presented approach operates on the simulation data only, providing a reduced surrogate model without the access to the original high-dimensional system matrices.

Constrained mechanical systems are very
often modeled with DAE systems that have an ODE representation in a reduced subspace on a constrained manifold. Operator inference provides the respective ODE system realization, as long as the solution data, used for the reduced-order model construction, stays on the constrained manifold. The algebraic constraints are eliminated by projecting the DAE snapshot data onto the POD subspace, constructed according to the constraints type. Moreover, the inferred reduced-order models retain the stability and interpretability properties due to the symmetric positive definite constraints imposed during the optimization process.

Additionally, we discussed the simulation parameters for the data generation that can affect the error when approximating the original system dynamics. We choose a non-harmonic excitation signal that allowed us to make an informed decision about the simulation duration and avoid a leap in error magnitude. However, characterizing the impact of the snapshot dataset on the quality of the reduced-order model remains an open question, also for the presented application case.
\section*{Acknowledgments}
The authors acknowledge the support and computational resources provided by the Max Planck Institute for Dynamics of Complex
Technical Systems, Magdeburg. This research was also supported by the Research Training Group ‘‘Mathematical Complexity
Reduction’’, which is a Graduiertenkolleg (DFG-GRK 2297) funded by Deutsche Forschungsgemeinschaft (DFG), Germany.
\addcontentsline{toc}{section}{References}
\bibliographystyle{plainurl}
\bibliography{journals,mor-adds}

\end{document}